\documentclass[11pt]{article}
\usepackage{latexsym}
\usepackage{graphicx}
\usepackage{amsfonts}
\usepackage{latexsym,amssymb,amsmath}
\usepackage[applemac]{inputenc}
\usepackage{graphicx}
    \usepackage{url}

\newcommand{\proscal}[2]{\langle #1,#2 \rangle}

 \newcommand{\HH}{\mathcal{H}(\Omega)}
 
  \newcommand{\ldeux}{L^2(\Omega)}

  \newcommand{\hun}{H^1(\Omega)}
   \newcommand{\demo}{\noindent\textit{ Proof -~}}
 \newcommand{\findemo}{\hfill $\Box$}
 \newcommand{\nd}{\noindent}
 \newcommand{\tr}{\ensuremath{\mathop{\rm Tr\,}\nolimits}}
 \newcommand{\Section}{\setcounter{equation}{0} \section}

\def\R{{\mathbb{R}}}
  
\newtheorem{theoreme}{Theorem}[section]
\newtheorem{definition}{Definition}[section]

\newtheorem{prop}{Proposition}[section]
\newtheorem{lemma}{Lemma}[section]

\newtheorem{rem}{Remark}[section]

\begin{document}
\title{Existence and uniqueness results for the Gradient Vector Flow and geodesic active contours mixed model}
\author{Laurence GUILLOT$^*$, Maïtine BERGOUNIOUX$^*$} 
\date{\today}
\maketitle
\begin{center}
\emph{$^*$UMR 6628-MAPMO, Fédération Denis Poisson, Université d'Orléans, BP. 6759, F-45067 Orléans Cedex 2,
laurence.guillot@univ-orleans.fr, maitine.bergounioux@univ-orleans.fr}
\end{center}


\begin{abstract}
This article deals with the so called GVF (Gradient Vector Flow) introduced by C. Xu, J.L. Prince \cite{Xu97, Xu98}. We give existence and uniqueness results for the front propagation flow for boundary extraction that was initiated by Paragios, Mellina-Gottardo et Ralmesh \cite{PMGR01,PMGR04}. The model combines the geodesic active contour flow   and the GVF to determine the geometric flow. The motion equation is considered within  a level set formulation to result an Hamilton-Jacobi equation. 
\end{abstract}


\Section{Introduction}
We consider an image  segmentation model : $I$ is  a given image and we want to detect   boundaries without  connexity or convexity  assumptions on  contours. Therefore we are interested in the Gradient Vector Flow (GVF) as a front propagation flow model. This model builds a class of vector fields derived from images and has been  introduced by Chenyang Xu and Jerry L. Prince in \cite{Xu97}. The GVF can be viewed as external forces for active contour models: it allows  to solve problems where classical methods convergence fail to deal with boundary concavities. 
On the other hand, a new front boundary-based geometric flow for boundary extraction was proposed by Paragios, Mellina-Gottardo and Ralmesh in \cite{PMGR01,PMGR04}. In this model the GVF is used to revise the geodesic active contour model of V. Caselles, R. Kimmel, G. Sapiro \cite{GAC} resulting on a bidirectional geometric flow.  
The classical parametric active contour model was proposed by D.Terzopoulos, A.Witkin et M.Kass \cite{Terzopoulos}. It derives  from    an energy functional minimization. Curves are drawn toward the object boundaries under  potential external forces  action which can be written as the negative gradient of a scalar potential function derived from images. Other forces (as pressure forces) may  be added. However, the performance of methods is limited by   unstable initialization  process and poor convergence when the boundary is concave.   Chenyang Xu and Jerry L. Prince \cite{Xu97} set a new external force, dealing with these limitations.

Most of time  snakes model provide a local minimum of the functional cost.  So,  C.Xu et JL.Prince  generalized this model from the balance equation at the equilibrium between internal and external forces (Euler equation of the minimization problem).  They replace the standard external force $F_{ext}^{(p)}=-\nabla P_{\bf{x}}$ (where $P$ is a potential edges detector), by a more general  external  force.
 
This new  external force field is  called the Gradient Vector Flow (GVF). It  includes a divergence-free component  and a curl-free component \cite{Kichen,Xu-Yezzi}. Therefore, this new active contour model cannot be formulated as an energy minimization problem.   The external force, denoted  $V$ below  is introduced  via the balance equation  that can be written :
\begin{equation}\label{eqn:GVFsnake1}
\alpha  \frac{d^2{\bf{x}}}{ds^2} - \beta\frac{d^4{\bf{x}}}{ds^4}+V({\bf{x}})=0
\end{equation}

The resulting parametrized curve solving the equation (\ref{eqn:GVFsnake1}) is called `` GVF-snake ''.  The final configuration of a GVF-snake satisfies an equilibrium equation which is not a variational problem Euler equation, since $V(x,y)$ is not an irrotational field.

In  next section, we introduce the GVF  and present the advantages of such a field in view of an  active contour model and we present and give a precise definition of the Gradient Vector Flow. 
In  section 3. we present this model which is inspired by the geodesic active contour model combined with the GVF. Finally, the level set strategy leads to  an Hamilton-Jacobi equation that has not been studied  yet (to our knowledge): in   the last section we give existence and uniqueness results.
\Section{The Gradient Vector Flow (GVF)}
The GVF $V$ is  a 2-dimensional vector field that should minimize the following objective function \cite{Xu97,PMGR01}.  
\begin{equation}\label{eqn:defGVF}
E(V) : =  \int_{\Omega} \mu(u_x^2+ u_y^2 +v_x^2+ v_y^2) + f|\nabla f|^2 |V- \nabla f|^2 d{\bf{x}} 
\end{equation}  

 \noindent where $\Omega$ is an open, bounded subset of $\R^2$, $ u_x, u_y,v_x,v_y $  denotes the spatial derivatives (with respect to $x$ and $y$)  of the field $V=(u,v)$, 
$\mu >0$ and $f : \Omega \to \R$ a continuous edge detector.

There are  many  choices for $f$:  
in \cite{Xu98}, C.Xu and JL.Prince consider
$$f(x,y)= E_{ext}^{(1)} (x,y) = -|\nabla I(x,y)|^2$$
where $I$ is the (intensity) image or 
$$f(x,y)=E_{ext}^{(2)} (x,y) = -|\nabla [G_{\sigma}*I(x,y)]|^2$$
where   $ G_\sigma $ is a Gaussian kernel, so that the image is filtered. 
Here we choose detector  proposed by R. Deriche and  O. Faugeras in \cite{Der-Faug}  
$$ f(x,y)=h(|\nabla(G_\sigma * I)(x,y)|^2) $$
where  \begin{equation} \label{defh}
 h(b)=1-\frac{1}{\sqrt{2\pi}\sigma }e^{-{\frac{b}{2\sigma ^2} }}~.\end{equation} 
Here  $|\cdot |$ denotes the  $\R^2$- euclidean norm. \\
The first term of the functional $E$ is a regularization term whereas the second term is a data-driven component. 
If $f|\nabla f|^2 $ is small, the energy is dominated by the first term and we get  a slowly varying field. On the hand, when $f|\nabla f|^2 $ is large, the second term forces  $V$ to decrease to $ \nabla f $.
Therefore $V$ is close to the gradient of the edge map when it is large (this is classical external force for snakes) and  does not evolve quickly  in homogeneous regions (which increase the  snake capture area).
The parameter $\mu$ govern the tradeoff between the two integrands of the cost functional and  should be set according to the image   noise level   (more noise increases $\mu$) see \cite{Xu98}.
\subsection{On GVF   existence }
 
 First we precise the notations :
 $\Omega$ is an open, bounded subset of $\R^2$ with $\mathcal{C}^\infty $ boundary  $\Gamma$: it is the image domain. 
 Let $V=(u,v)$ and  $W=(\xi,\chi)$  in  $\HH:= H^1(\Omega) \times H^1(\Omega)$; then   $ DV{\bf{(x)}}= (  \nabla u{\bf{(x)}} ,  \nabla v{\bf{(x)}} )  \in {\ldeux}^4$. 
The inner product in ${\ldeux}^4$ is
 $$\proscal{DV}{DW}_2= \int_\Omega \proscal{\nabla u}{\nabla \xi}_{\R^2}+\proscal{\nabla v}{\nabla \chi}_{\R^2} d{\bf{x}}  =\proscal{\nabla u}{\nabla \xi}_{{\ldeux}^2}+\proscal{\nabla v}{\nabla \chi}_{{\ldeux}^2},$$
and  the ${\ldeux}^4$-norm of $DV$    is denoted   $ \|DV\|_2 $.
 The space $\mathcal{H}(\Omega)= H^1(\Omega) \times H^1(\Omega)$ is endowed with the  norm   $$\|V\|^2_{\mathcal{H}(\Omega)}=\|V\|^2_{{\ldeux}^2}+\|DV\|^2_2.$$

 A first  ``definition'' of the GVF could be the following:
{\definition{
The Gradient Vector Flow field is defined as a solution of the following optimization problem 
\begin{equation} \label{eqn:PbMinim}
 \mathcal{(P)} :  \qquad  \min \{ E(V)~|~  V \in \HH \}. 
 \end{equation}
where the energy functional $E(V)$ is defined by  (\ref{eqn:defGVF}).}}   
Unfortunately,  such a definition is not correct since the minimum is not necessarily attained. The study of $E$  will provide another  definition of the GVF as the  solution to  a decoupled system parabolic partial differential equations.  

From now and in the sequel we assume  that the edge function $f$ verifies  $(\mathcal{H}_1)$  :
$$(\mathcal{H}_1) : \qquad  f \in \mathcal{C}^{1}( \overline{\Omega}) \mbox{ and } 
f \geq 0$$
  
  \begin{rem} Note that if $f$ satisfies  $(\mathcal{H}_1)$ then  $ f\in H^1(\Omega) $  and  $ f |\nabla f|^2 \in L^{\infty}(\Omega)$ which are the  ``minimal'' assumptions in a first step.
  \end{rem}
  
{\prop : Assume $(\mathcal{H}_1)$. The functional $E$ is  continuous on $\HH$.}
\\
\demo We get   $$E(V)=\mu\|DV\|_2^2+\int_\Omega f|\nabla f|^2(|V|^2-2\proscal{V}{\nabla f}_{\R^2}) d{\bf{x}}+\int_\Omega f|\nabla f|^4 d{\bf{x}} =\pi(V,V) + L(V)+C  $$
 where $$\pi(X,Y)=\mu \proscal{DX}{DY}_2+\int_\Omega f|\nabla f|^2\proscal{X}{Y}_{\R^2} d{\bf{x}}, \hspace{0.3cm} (X,Y)\in (\HH)^2~,$$  $$L(X)=  - 2\int_\Omega f|\nabla f|^2 \proscal{X}{\nabla f} d{\bf{x}}, \hspace{0.3cm}X\in\HH \mbox{ and } C =\int_\Omega f|\nabla f|^4 d{\bf{x}}~.$$
The bilinear form $\pi$ is continuous on $\HH$ for the $ \hun^2$-norm : 
let $X$ et $Y$ be  in $\HH$,
\begin{align}
  |\pi(X,Y)|&=  \left|\mu \proscal{DX}{DY}_2+\int_\Omega f|\nabla f|^2\proscal{X}{Y} d{\bf{x}}\right| \notag\\
 &\leq  \mu \|DX\|_2 \|DY\|_2+\|f|\nabla f |^2\|_{L^\infty(\Omega)}  \|X\|_{{\ldeux}^2} \|Y\|_{{\ldeux}^2} \notag\\
  & \leq   C(\mu,f) \|X\|_{{\hun}^2}\|Y\|_{{\hun} ^2} \notag\end{align}
where $C(\mu,f)$ is a constant that only depends on data.
 $L$ is  obviously linear and continuous on $ \HH $. We deduce the continuity of $E$.   \findemo

 {\theoreme
The functional $E$   is Gâteaux-differentiable on $\mathcal{H}(\Omega)$ and   for every $ V=(u,v)$ and $W=(\xi,\chi) $ in $\HH $  we get 
\begin{equation} \label{diffE}
\proscal{\nabla E(V)}{W}=2  \int_\Omega  \mu \proscal{DV}{DW}+ f|\nabla f|^2  \proscal{V-\nabla f}{W}d{\bf{x}} \end{equation}
Furthermore, if $\mu\in \R^+$ then $E$ is convex on $\HH$.
}
\\
\demo The Gâteaux -differentiability is clear.  In addition for every 
 $V=(u,v)$ and  $W=(\xi,\chi)$  in $\HH$ we have
 $$  
E(V+W)-E(V)-\proscal{\nabla E(V)}{W}= \int _\Omega \mu |DW|^2 + f|\nabla f|^2|W|^2 d{\bf{x}}~,
$$
so that,  if $\mu\geq 0$,  as $f\geq 0$ by $(\mathcal{H}_1)$
 $E(V+W)-E(V)-\proscal{\nabla E(V)}{W} \geq 0.$ So we get the convexity of $E$  on $\HH$. 
\findemo
 
 ~
 
\noindent\textbf{From now we assume $\mu> 0.$ }
{\theoreme
Assume that  $f|\nabla f|^2$ does not  degenerate on $\Omega$ : 
\begin{equation} \label{hplus} \exists c >0~~ \hbox{ such as } f|\nabla f|^2 \geq c  \end{equation}
then  $\pi$ is coercive on $\HH$. This  implies the coercivity and  the strict-convexity of  $E$ on $\HH $ and  the problem $(\mathcal{P})$  has  a unique solution.
}
\\
\demo  $\pi$ has been define above and we get 
$$\pi(X,X) = \mu \|DX\|_2^2+\int_{\Omega}f|\nabla f|^2|X|^2_{\R^2} d{\bf{x}}    
 \geq \mu\|DX\|_2^2+c \int_{\Omega}|X|^2_{\R^2} d{\bf{x}}  
 \geq \min(\mu,c)\|X\|_{{\hun}^2}~; 
$$
therefore $\pi$  is ${{\hun}^2}$-coercive on $\HH$ and strictly convex. Therefore  $E$ is coercive and strictly convex on $\HH$ as well.
It follows that if $(\ref{hplus})$ is verified then $(\mathcal{P})$ admits a unique solution. 
 \findemo

 Let us formally write the Euler-equations  of problem $ (\mathcal{P}) $: assume 
 $V^*=(u^*,v^*)$ is a solution to  $( \mathcal{P} )$. By convexity,  it is   a stationary point and    $\nabla E(V^*)=0$. 
 Let $W=(\xi,\chi) \in \mathcal{H}(\Omega).$ 
Integrating by parts  expression (\ref{diffE} )  gives
\begin{align}\label{eqn:numero1}
&\proscal{\nabla E(V^*)}{W} \notag \\&=- \mu  \int_\Omega(\xi \Delta u^*   +\chi\Delta v^*   ) d{\bf{x}}+ \int_\Omega  f|\nabla f|^2  \proscal{V^*-\nabla f}{W}_{\R^2}d{\bf{x}} + \int_\Gamma \left(\frac{\partial u^*}{\partial \nu} \xi +\frac{\partial v^*}{\partial \nu} \chi \right) d\sigma
\end{align}
where $\nu({\bf{x}})$ is the outer unit normal of $\partial \Omega=\Gamma$ at ${\bf{x}}$.
\\
We  first suppose that $\xi$  and  $\chi$  belong to  $\mathcal{D}(\Omega)$ (the space of $C^\infty (\Omega)$ functions with compact support in $\Omega$). Then
 $$- \mu  \int_\Omega(\xi \Delta u^*   +\chi \Delta v^*   ) d{\bf{x}}+ \int_\Omega  f|\nabla f|^2  \proscal{V^*-\nabla f}{W}_{\R^2}\, d{\bf{x}}=0 $$
So  
$$\left \{\begin{array}{ll}
- \mu \Delta u^*+f|\nabla f|^2(u^*-f_x)=0 \hbox{ in } \mathcal{D'}(\Omega) \\[0.3cm]  -\mu \Delta v^* + f|\nabla f|^2(v^*-f_y)=0 \hbox{ in } \mathcal{D'}(\Omega) \end{array}\right.$$
\indent where $f_x$, $f_y$ are spatial derivatives of $f$.
With  (\ref{eqn:numero1}), we obtain : 
 $$\proscal {\left( \begin{array}{ccc}\displaystyle{\frac{\partial u^*}{\partial \nu}} \\ \displaystyle{\frac{\partial v^*}{\partial \nu}} \end{array} \right)}{\left( \begin{array}{ccc} \displaystyle{\xi |_{\Gamma}}\\ \displaystyle{\chi |_{\Gamma} }\end{array} \right)}=0$$
Finally, if  a solution $V^*=(u^*,v^*)$ to    problem $(\mathcal{P})$  exists it must verify : 
\begin{equation} \left \{ \begin{array}{c}
\label{eqn:uellip} -\mu \Delta u^*+f|\nabla f|^2(u^*-f_x) =0~\mbox{ in } \Omega , \displaystyle{\frac{\partial u^*}{\partial \nu}}=0 
\mbox{ on } \Gamma \\[0.3cm]
-\mu \Delta v^* + f|\nabla f|^2(v^*-f_y)=0 ~\mbox{ in } \Omega , \displaystyle{\frac{\partial v^*}{\partial \nu}}=0  \mbox{ on } \Gamma ~. \end{array}\right.
\end{equation}
 
The above equations are equilibrium equations  if $V^*$ realizes the minimum of  energy $E$.  However, we cannot ensure the existence of such a minimum. Indeed, assumption (\ref{hplus}) is not realistic :  the same gray level for  image $I$ on a significant area  implies $\nabla f=0$.
In this case   the functional $E$ is  a priori non longer  coercive and   we do not know if $ (\mathcal{P})$ has  a solution. So, instead  of computing the ``exact''   GVF, minimizer of the functional $E$, we   approach it by a minimizing sequence and we consider  it is the equilibrium state of a  time evolving vectors fields. The stationary problem becomes a  {\bf{ dynamic}} one : 
 $$V^*({\bf{x}})=\lim_{t\to +\infty}V(t,{\bf{x}}).$$
 A simple way to impose a motion to the vectors field is to impose the velocity $\displaystyle \frac{\partial V}{\partial t}({\bf{x}},t)$  setting :  
\begin{equation}\label{eqn:vitesse}
\frac{\partial V}{\partial t}({\bf{x}},t)=-\nabla_V\, E(V(t,{\bf{x}}))
\end{equation}


This leads to  parabolic partial differential equations. The gradient vector flow is initialized as the gradient of the edge detector $f$ :
\begin{equation}
V(t=0,{\bf{x}})=V_0({\bf{x}})=\nabla f({\bf{x}}) \hbox{ sur } \Omega 
\end{equation}
 We obtain the  following {\bf{dynamic formulation }} which  is the GVF suitable  definition : 
 {\definition{The gradient vector flow  $V=(u,v)$ is defined by the following decoupled equations respectively verified by each of its coordinates $u$ and $v$ : 
\begin{equation} \label{eqn:u}
\left \{ \begin{array}{ll}
 \displaystyle{  \frac{\partial u}{\partial t}} - \mu \Delta u -f|\nabla f|^2  (u-f_x)  = 0 &\mbox{in } ]0, +\infty [\times \Omega \\
 \displaystyle{\frac{\partial u}{\partial \nu}}=0 &\mbox{on } ]0, +\infty [\times \Gamma\\
 u(0,\cdot)=  f_x &\mbox{in }\Omega \end{array}\right.\end{equation}
 \begin{equation} \label{eqn:v}
 \left \{ \begin{array}{ll}
 \displaystyle{  \frac{\partial v}{\partial t}} - \mu \Delta v -f|\nabla f|^2  (v -f_y)  = 0 &\mbox{in } ]0, +\infty [\times \Omega \\
 \displaystyle{\frac{\partial v}{\partial \nu}}=0 &\mbox{on } ]0, +\infty [\times \Gamma\\
 v(0,\cdot)= f_y &\mbox{in }\Omega \end{array}\right.\end{equation}
 }}
  
\subsection{GVF properties.}\label{GVF properties}

In this subsection,  we give regularity  properties of solutions to   (\ref{eqn:u}) and  (\ref{eqn:v}). Of course, it is sufficient to consider equation  (\ref{eqn:u}).  First, the existence of a unique solution is given by a classical theorem (see for example  \cite{LM3,Fr1,Fr2}).
The  bilinear form $a$ associated to equation  (\ref{eqn:u}) is 
$$a(t;u,v)=\mu \int_\Omega \left(\nabla u({\bf{x}}).\nabla v({\bf{x}}) +u({\bf{x}}) v({\bf{x}})\right) d{\bf{x}}+\int_\Omega f|\nabla f|^2 u({\bf{x}}) v({\bf{x}})d{\bf{x}}$$ where $\mu>0$ and $ f|\nabla f|^2\in L^{\infty}( ]0,T[ \times \Omega)$.
It satisfies : 
\begin{enumerate}
\item  $t \mapsto a(t;u,v)$ is measurable $ \forall u,v \in H^1(\Omega)$,
\item For almost $t \in[0,T]$ and for all  $u,v \in H^1(\Omega)$,
$$|a(t;u,v)|\leq C(\mu,f) \|u\|_{H^1(\Omega)}\|v\|_{H^1(\Omega)}$$
 \item For almost  $t \in [0,T]$ and for all $v \in H^1(\Omega)$,
 $$a(t;v,v)\geq \mu \|v\|_{H^1(\Omega)}^2 +\|f|\nabla f|^2\|_{L^\infty(\Omega)} |v|_{L^2(\Omega)}^2$$
  \end{enumerate}
Note that here $a$ does not depend on $t$. So we may assert that  (\ref{eqn:u})  has a unique solution 
 $u \in L^2(0,T;H^1(\Omega))\cap \mathcal{C}([0,T];\ldeux) $  and $\displaystyle{\frac{\partial u}{\partial t} \in L^2(0,T; H^1(\Omega)').}$ Moreover  
 
{\theoreme \label{regulGVF}{
Let $T >0$, and assume ($\mathcal{H}_1).$  Then the GVF  $(u,v)$ (solution of   (\ref{eqn:u}) and  (\ref{eqn:v}))
 is $\mathcal{C}^1$ on $[0,T] \times \overline{\Omega}$.
}}\\
\demo We use a generic regularity result (\cite{LM3}). We prove the  result for the component $u$. Assumption $(\mathcal{H}_1)$ yields that $\nabla f \in L^2(\Omega)^2$ and $ f |\nabla f|^2 \nabla f \in [\mathcal{C}^1([0,T] \times \overline{\Omega})]^2$. So the solution $u$ of  ( (\ref{eqn:u}) 
belongs to $\mathcal{C}^ 1([\varepsilon,T] \times \overline{\Omega} )$ for all $\varepsilon>0$.
Moreover $\nabla f \in \mathcal{C}^0 (\overline{\Omega})$ according to $(\mathcal{H}_1)$ and  compatibility  conditions are satisfied with respect to boundary and initial data. So we may conclude. \findemo

The GVF is built as a spatial diffusion of the  edge detector $f$ gradient.
This  is equivalent to a progressive construction of the gradient vector flow starting from the object boundaries and moving toward the flat background. In \cite{Xu98}, the GVF is normalized to obtain a more efficient propagation. It is denoted 
 $\hat V{\bf(x)}=(\hat u {\bf(x)},\hat v {\bf(x)})$ where 
$$ 
 \hat u {\bf(x)} = \frac{u{\bf(x)}}{\sqrt {u{\bf(x)}^2+v{\bf(x)}^2}} , \hspace{ 0.7cm } \hat v {\bf(x)}= \frac{v{\bf(x)}}{\sqrt {u{\bf(x)}^2+v{\bf(x)}^2}}
$$ 
 to give  the new external force of the geometric flow called \textbf{GVF-snake mode}l. The velocity of the contour $C$   is given by the equation : 
\begin{equation}\label{eqn:GVFsnake}
C_t(s,t)=\alpha  \frac{d^2{\bf{x}}}{ds^2} - \beta\frac{d^4{\bf{x}}}{ds^4}+\hat V({\bf{x}})
\end{equation}
C. Xu and J.L. Prince have shown that a such flow is  not dependent of   initial conditions and deal with concave regions. However, it depends on   curve parametrization, cannot manage topological changes, and involves second and fourth order derivatives that are difficult to estimate. The objective of N. Paragios, O. Mellina-Gottardo and  V. Ralmesh in \cite{PMGR01} is to eliminate these shortcomings by integrating the GVF with the geodesic active contour and implementing it using the level set method \cite{PMGR04}. Our aim is to study the Hamilton-Jacobi equation derived from this model.
\section{Geodesic contours and GVF} 
\subsection{ Paragios - Mellina-Gottardo - Ralmesh model  }
N. Paragios, O. Mellina-Gottardo and  V. Ralmesh have defined in \cite{PMGR01} a new ``front propagation flow for boundary extraction''. Their geometric model is inspired by the geodesic active contours (\cite{GAC}) and directly defined by the contour evolution velocity. It is based on the remark that the Gradient Vector Flow field after the rescaling refers to the direction that has to be followed to locally deform the contour  and to reach  the closest object boundaries.  On the other hand, given the fact that the propagation of a contour often occurs along the normal direction,  the propagation will be optimal when $\hat{V}$ and the unit outward normal $\nu$ are colinear. So  we choose to project the normalized gradient vector flow onto the outward normal. Then we multiply the velocity by an edge detector function $g$ (that may be different from $f$), which represents the contour information. The contour evolution velocity is then given by the equation : 
\begin{equation} \label{eqn:GVF-geod1}
\mathcal{C}_t({\bf{x}}) =\underbrace {g(|\nabla I_\sigma({\bf{x}})|^2) }_ {\mbox{\tiny{boundary }}} \underbrace {\proscal{(\hat{u},\hat{v})({\bf{x}})}{ \nu({\bf{x}})}_{\R^2}}_ {\mbox{\tiny{projection}}} \nu({\bf{x}})
\end{equation}
 where $I_\sigma$ is the filtered image and $ g(b)=\frac{1}{\sqrt{2\pi}\sigma}e^{-{\frac{b}{2\sigma^2} }} = 1-h(b) $ where $h$ is defined by (\ref{defh})
\\
When there is no  boundary information ( $|\nabla I_\sigma|^2  \ll 1$), the contour evolution is driven by the inner product between the Normalized Gradient Vector Flow (NGVF) and the normal  direction : it is adapted to deal with concave regions. When  the curve reaches   the object boundaries  neighbourhood ( $|\nabla I_\sigma|^2 \simeq + \infty $ ) then    $g\approx 0$ that is, the flow becomes inactive and the equilibrium state is reached.

It is classical to impose a regularity condition on the contour propagation  adding  a curvature term and a ``balloon force'' $H$. The evolution equation becomes  
 \begin{equation}
 \mathcal{C}_t({\bf{x}}) =g(|\nabla I({\bf{x}})|^2)\left(\underbrace {-\beta \kappa ({\bf{x}})}_{\mbox{\tiny{ smoothness}}} +\underbrace { (1- |H({\bf{x}})|)  \proscal{(\hat{u},\hat{v})({\bf{x}})}{\nu({\bf{x}})}_{\R^2} } _{\mbox{\tiny{ boundaries attraction }}}
 +\underbrace { H({\bf{x}})}_{\mbox{\tiny{ balloon force}}}\right)\nu({\bf{x}})
\end{equation}
where $\beta >0$. 
 
\subsection{Level set implementation}
Problems of  topologic changes can be solved using the level set method \cite{OS}. The moving 2D- curve  is viewed as the zero  level set of a 3D surface which  equation is $z-\Phi(x,y)=0$. We denote $\Phi_t$ the partial derivative $\displaystyle{\frac{\partial \Phi}{\partial t}}$ of $\Phi$ towards $t$.
 
{\theoreme
The evolution of the  3D-surface $\Phi$ is described by : 
\begin{equation} \label{eqn:levelsetparagios}
\left\{\begin{array}{lll}
 \Phi_t({\bf{x}}) &= g(|\nabla I_\sigma({\bf{x}})|^2) \left( (\beta \kappa ({\bf{x}})-H({\bf{x}}))|\nabla \Phi ({\bf{x}})|\right.& \\
  &~~~~~~\left.- (1- |H({\bf{x}})|)\proscal{ \hat V ({\bf{x}})}{\nabla \Phi({\bf{x}})}_{\R^2} \right) &  \mbox{in }  ]0,+\infty[ \times \Omega  , \\
 \displaystyle{\frac{\partial \Phi}{\partial \nu}} &=0   & \mbox{on } ]0,+\infty[ \times \Gamma ,\\[0.3cm]
  \Phi(0,\cdot)&=\Phi_o   &\mbox{in }\Omega 
\end{array} \right.
\end{equation}
where  $\Phi_o$ is the\textbf{ signed distance}   defined by : 
$$  \Phi_o(\bf{x}) = ~\left\{ \begin{array}{l} 
 0 \mbox{ if }{ \bf{x} } \in \Gamma ~,\\[0.2cm]
  \pm ~d  \mbox{  the distance between } { \bf{x}}   \mbox{ and } \Gamma~; 
 \end{array} \right.
$$
the positive sign (resp. negative) is chosen if the point $\bf{x}$ is outside (resp. inside)  $\Gamma $.}
\\
\demo  The curve $\mathcal{C}$ is the zero set level of $\Phi$ : 
 $  \Phi( t,\mathcal{C}({\bf{x}})) = 0 .$ 
 Deriving  formally with respect to  $ t$  gives  for almost every ${\bf{x}}\in \Omega$:
$$
 \Phi_t ({\bf{x}})+ \proscal{\nabla \Phi({\bf{x}})}{\mathcal{C}_t ({\bf{x}})}_{\R^2}= 0
$$
$$
\Phi_t ({\bf{x}})+ \proscal{\nabla \Phi ({\bf{x}})}{g(|\nabla I_\sigma({\bf{x}})|^2)\left(-\beta \kappa ({\bf{x}}) +(1- |H({\bf{x}})|)  \proscal{\hat V({\bf{x}})}{\nu({\bf{x}})}_{\R^2}+H({\bf{x}}) \right)\nu({\bf{x}})}_{\R^2} = 0
$$
where $ \nu({\bf{x}}) = \displaystyle{\frac{\nabla\Phi({\bf{x}})}{|\nabla\Phi({\bf{x}})|}} $ is the outward  unit normal. 
A short computation gives 
 $$
\Phi_t ({\bf{x}})= g(|\nabla I_\sigma({\bf{x}})|^2)\left(\beta \kappa ({\bf{x}})- (1 - |H({\bf{x}})| ) \proscal{\hat V({\bf{x}})}{\displaystyle{\frac{\nabla \Phi({\bf{x}})}{|\nabla \Phi({\bf{x}})|}} } -H({\bf{x}})\right)  |\nabla \Phi({\bf{x}})| 
$$
Finally  \vspace{-0.8cm}
$$
\Phi_t({\bf{x}}) = g(|\nabla I_\sigma({\bf{x}})|^2) \left( (\underbrace{\beta \kappa ({\bf{x}})}_{\mbox{(a)}}-\underbrace{H({\bf{x}})}_{\mbox{(c)}})|\nabla \Phi ({\bf{x}})| - (1- |H({\bf{x}})|)\underbrace{\proscal{ \hat V ({\bf{x}})}{\nabla \Phi({\bf{x}})}}_{\mbox{(b)}}\right) 
$$ 
 \findemo

 The final flow can be decomposed in 
 \\
$\bullet  $ (a) :  a term that provides propagation  regularity   aims  and shrinks  the curve toward the object boundaries,
\\
$\bullet  $ (b) :  a bidirectional flow that moves the curve toward the internal and external objects boundaries,
\\
$\bullet  $ (c) :  an adaptative balloon force that drives the propagation of the curve when the boundary term becomes inactive \cite{PMGR04}.

\subsection{Choice of the edge detector $f$  }

As mentionned before, we  focus on the Gaussian edge detector proposed by R. Deriche et O. Faugeras in \cite{Der-Faug}: 
  \begin{equation}\label{eqn:deff}
  f(x,y)=h(|\nabla(G_\sigma * I)(x,y)|^2) \end{equation}
  where  $h$ is defined by (\ref{defh}).
 We must  ensure that  assumption $(\mathcal{H}_1)$ is satisfied for the edge detector $f$. Therefore, we have to set  additional  hypothesis   the image  intensity function that we have denoted  $I$.

$I$ is supposed to have  compact support included in   $\Omega$  that is   image  frame:  we decide that $I$ is equal $0$ out of the ``true image".
$I\in L^{\infty}(\Omega)\subset L^{2}(\Omega)\subset L^{1}(\Omega).$
On the other hand $I\notin H^1(\Omega)$. Indeed, the image gradient norm  becomes  infinite  at objects contours. Furthermore $I\notin \mathcal{C}^{1}(\Omega)$.  We correct this  lack of regularity using a filtering process that makes the filtered image $\mathcal{C}^{\infty}$ (and of course $\mathcal{C}^{1}$).

The filtered image   is  supposed  to be $\mathcal{C}^{\infty}$, with compact support.  So, we cannot choose  $G_\sigma * I$ since the resulting image has no  compact support.  So we consider a fixed  compact subset   $X$ of $\Omega$ and $Y$ a  compact subset  of $\Omega$ containing $X$. \\
  Let us consider $G^X_{\sigma}$ a  $\mathcal{C}^{\infty}(\Omega)$  projected of the Gaussian kernel $G_\sigma$ such as :  
$$G^X_{\sigma}({\bf{x}})=\left \{ \begin{array}{lc}
G_{\sigma} ({\bf{x}}) \hbox{ if } {\bf{x}}\in X \\  
0  \hbox{      if }{\bf{x}}\notin Y
\end{array}\right.$$
\nd So the regularized (filtered) image $I_{\sigma}=G_{\sigma}^X\ast I$ verifies $$I_\sigma \in \mathcal{C}^\infty_c(\Omega)$$ and $$\forall {\bf{x}} \in \Omega, \nabla I_\sigma({\bf{x}})=\nabla (G_\sigma^X \ast I)({\bf{x}})=\nabla G_\sigma^X  \ast I ({\bf{x}}).$$
\nd Of course the support of $I_{\sigma}$, contained in $\overline{\textrm{supp}(G_{\sigma}^X) + \textrm{supp}(I)}$, is not necessarily included in $\Omega$, but even if we must extend the frame $\Omega$ of the image, we consider that the filtered image has a compact support in $\Omega$. More precisely : {\prop \label{eqn:gradI}
The function 
\begin{eqnarray}
\left\{\begin{array}{ll}
\Omega  \to \R \\
{\bf{x}}\mapsto \nabla I_{\sigma}{\bf{(x)}}
\end{array} \right.
\end{eqnarray}
is $\mathcal{C}^\infty(\Omega)$, with compact support in $\Omega$
and so bounded on $\Omega$.
$$\nabla I_{\sigma} \in L^{\infty}(\Omega)\cap\mathcal{C}^\infty_c(\Omega)$$
}
From now and in the sequel, we denote $I_{\sigma}$ by $I$ and we make the following hypothesis on $I$: 
$$(\mathcal{H}_I):\quad I\in\mathcal{C}_c^{\infty}(\Omega) \cap W^{1,\infty}(\Omega)$$
Note that $I$ may be extended by $0$ to $\overline{\Omega}$.
Now we can give $f$ regularity properties : 
{\prop \label{regulf} Assume the image $I$ satisfies  $(\mathcal{H}_I)$.
The function $f$ defined by : 
\begin{eqnarray}
f: \left\{\begin{array}{ll}
\overline{\Omega} \to \R \\
{\bf{x}}\mapsto f({\bf{x}})=h(|\nabla I ({\bf{x}})|^2)=1-\displaystyle{\frac{1}{\sqrt{2\pi}\sigma}}e^{-{{\frac{|\nabla I ({\bf{x}})|^2}{2\sigma^2}}}}
\end{array} \right.
\end{eqnarray}
belongs to  $\mathcal{C}^{\infty}(\overline \Omega)$.
 } \\
\demo It is clear since  $I$  is $\mathcal{C}^\infty(\overline{\Omega})$, with compact support in $\Omega$ and $h$ is obvioulsy  $\mathcal{C}^\infty$. In addition 
    $f$ is   nonnegative and bounded by  $1$ on $\Omega$.
 \findemo 
 
Let us denote $\widetilde{g}$ the function  
$${\bf{x}}\in \overline{\Omega} \mapsto g(|\nabla I({\bf{x}})|^2)=\frac{1}{\sqrt{2\pi}\sigma}e^{-{\frac{|\nabla(I )({\bf{x}})|^2}{2\sigma^2} }}~.$$
We have seen in  Proposition  \ref {regulf} that   $\widetilde{g}=1-f  \in \mathcal{C}^\infty(\overline{\Omega})$.

In the sequel, we need the following lemma : 
{\lemma {\label {eqn:l2}  
The function $\sqrt{\widetilde{g}}$ belongs to $ \mathcal{C}^{1}(\overline{\Omega}) $ : it is Lipschitz continuous on $\overline{\Omega} $ with constant $K_1$.
}}\\ 

\section{Propagation equation study }
Equation (\ref{eqn:levelsetparagios})  has been obtained quite formally with the level-set method. Now, we give 
 existence and uniqueness  of a viscosity solution (\cite{GB} for example) .  Of course, we assume that the image $I$ verifies   $(\mathcal{H}_I)$ so that $f$  satisfies $(\mathcal{H}_1)$. The outward unit normal $\nu$  is  a  $ \mathcal{C}^{1,1}$ Lipschitz  vector field. 
 \subsection{Viscosity theory framework}
Let us precise the notations: 
for $p=\left(\begin{array}{c}p_1 \\ p_2\end{array}\right) \in \R^2$, the matrix $\displaystyle{\frac{1}{|p|^2}}\left(\begin{array}{cc}p_1^2 & p_1p_2 \\p_1p_2 & p_2^2\end{array}\right)$ is denoted by $ \displaystyle{\frac{p\otimes p}{|p|^2}}$. 
\\ 
We choose  $\beta\geq0$ and assume  for simplicity 
  that the curvature term $\kappa$ is the standard mean curvature  : 
 \begin{equation}\label{eqn:courbure}
 \begin{array}{c}
\displaystyle{\kappa({\bf{x}}) =\mbox{div }\left(\frac{\nabla \Phi}{|\nabla \Phi|}\right)({\bf{x}})}\\ [0.3cm] \displaystyle{=\left( \frac{\frac{\partial^2 \Phi}{\partial x^2}+\frac{\partial^2 \Phi}{\partial y^2}}{|\nabla \Phi|}-|\nabla \Phi|^{-1}\frac{\left(2\frac{\partial \Phi}{\partial x}\frac{\partial \Phi}{\partial y}\frac{\partial^2\Phi}{\partial x \partial y}+\left(\frac{\partial \Phi}{\partial y}\right)^2 \frac{\partial^2 \Phi}{\partial y^2}+\left(\frac{\partial \Phi}{\partial x}\right)^2 \frac{\partial^2 \Phi}{\partial x^2}\right)}{|\nabla \Phi|^2} \right)({\bf{x}})
}\\ \displaystyle{ =  \frac{\tr (D^2\Phi({\bf{x}}))}{|\nabla \Phi|} - |\nabla \Phi|^{-1} \tr \left( \frac{\nabla \Phi ({\bf{x}})\otimes \nabla \Phi ({\bf{x}})}{|\nabla \Phi ({\bf{x}})|^2}D^2\Phi({\bf{x}})\right)}
  \end{array}
\end{equation}
where $D^2\Phi$ is the Hessian matrix of $\Phi$.
Equation (\ref{eqn:levelsetparagios})   becomes 
  \begin{equation} \begin{split}
\Phi_t(t,{\bf{x}})&+g(|\nabla I({\bf{x}})|^2)  H({\bf{x}}) |\nabla \Phi ({\bf{x}})|
\\ &-g(|\nabla I({\bf{x}})|^2) \beta \left(\tr (D^2\Phi({\bf{x}}))-\tr \left( \frac{\nabla \Phi ({\bf{x}})\otimes \nabla \Phi ({\bf{x}})}{|\nabla \Phi ({\bf{x}})|^2}D^2\Phi({\bf{x}})\right)\right)
\\ &+g(|\nabla I({\bf{x}})|^2) (1- |H({\bf{x}})|)\proscal{\hat{V}({\bf{x}})}{\nabla \Phi({\bf{x}})}=0
\end{split}
\end{equation}
 So
 \begin{equation}\label{eqn:ls3}
\begin{split}
\Phi_t(t,{\bf{x}})&+ g(|\nabla I({\bf{x}})|^2)  H({\bf{x}}) |\nabla \Phi ({\bf{x}})|
\\ &-\beta g(|\nabla I({\bf{x}})|^2)  \tr \left( \left[ I- \frac{\nabla \Phi ({\bf{x}})\otimes \nabla \Phi ({\bf{x}})}{|\nabla \Phi ({\bf{x}})|^2} \right] D^2\Phi({\bf{x}}) \right)
\\ &+g(|\nabla I({\bf{x}})|^2)  (1- |H({\bf{x}})|)\proscal{\hat{V}({\bf{x}})}{\nabla \Phi({\bf{x}})}=0
\end{split}
\end{equation}
where $I$ is the identity matrix.
Setting
\begin{equation}\label {eqn:matA}
A(p):=I- \frac{p\otimes p}{|p|^2}.
\end{equation} 
and $F({\bf{x}},p,X) := $
\begin{equation}\label{eqn:F}
 g(|\nabla I({\bf{x}})|^2) H({\bf{x}}) |p|
 -\beta g(|\nabla I({\bf{x}})|^2)  \tr(A(p)X)
 +g(|\nabla I({\bf{x}})|^2) (1-|H({\bf{x}})|)\proscal{\hat{V}({\bf{x}})}{p}
\end{equation}
we get 
\begin{equation}\label{eqn:ls4}
\Phi_t+F({\bf{x}},\nabla \Phi ,D^2\Phi)=0 ~.
\end{equation}
The Hamiltonian $F:\R^2\times \R^2\times \mathcal{S}_2\to \R$, is  independent  of $t$ and $\Phi$. Here   $\mathcal{S}_N$ is  $N\times N$ the symmetric matrices space endowed with th classical  order. 
 \\
 First, we have to verify we may  use  the viscosity theory framework. 
 {\definition  Let $\mathcal{F}$ be a function : $ \left \{ \begin{array}{ccc}\R\times  \R^2\times   \R^2\times \mathcal{S}_2 &\to &\R  \\
 (t,x,p,X) & \mapsto & \mathcal{F}(t,x,p,X)\end{array}\right.$. \\ $\mathcal{F}$ is proper if :
 \\
 -  it is a non-decreasing function with respect to the first   variable  $t$ : \begin{equation}\label{eqn:propre2}
\mathcal{F}(t,x,p,X) \leq \mathcal{F}(s,x,p,X) \mbox{  if  } t\leq s. 
\end{equation}
- it is  a  decreasing function with respect to the  last variable $X$ 
\begin{equation}\label{eqn:ellip}
\mathcal{F}(t,x,p,X) \leq \mathcal{F}(t, x,p,Y) \mbox{  si  } Y\leq X. 
\end{equation}   
}
{\theoreme The function $ {F}$ is proper}
\\
\demo 
  Since $F$ does not explicitly depend on $t$, $F(t, {\bf{x}},p,X)=F({\bf{x}},p,X)$, the condition (\ref{eqn:propre2}) is verified.
\\ Assuming $p\neq0$. The matrix $A(p)$ defined by (\ref{eqn:matA}) is  semi-definite positive, and  can be written as $A=\sigma \sigma^t$.
So we have $$\tr(AX)=\tr(\sigma \sigma^t X)=\tr(\sigma^t X \sigma)=\sum_{i=1}^2 \sigma_i^t X \sigma_i,$$
where $\sigma_i$ stands for  the ith column of $\sigma.$\\
 Let $X,Y \in \mathcal{S}_2$ be given. Let us assume $Y\geq X$, then : 
 $$\forall i\in\{1,2\}, \sigma_i^t X \sigma_i \leq \sigma_i^t Y \sigma_i.$$
The constant $\beta$ is positive and the  function $g$ as well, so we deduce  that : 
 $$-\beta g(|\nabla I({\bf{x}})|^2) \tr(A(p)X)\geq-\beta g(|\nabla I({\bf{x}})|^2) \tr(A(p)Y).$$
 so $F({\bf{x}},p,X)\leq F({\bf{x}},p,Y)$. \\
We conclude that the function $F$ is  proper. 
\findemo 

\nd The  general viscosity theory  framework is thus well posed.

\subsection{Ishii and Sato theorem}
 
 In what follows we use a theorem by   Ishii and Sato \cite{IS} that gives existence of viscosity solution to singular degenerate parabolic  (Hamilton-Jacobi) equations with nonlinear oblique derivative boundary conditions.

In this subsection we recall this theorem.  In the sequel we denote 
  $\rho$ the function defined from $\R^2$ in $\R$ by $\rho(p,q)=\min\left(\displaystyle{\frac{|p-q|}{\min(|p|,|q|)},1}\right).$

 {\theoreme  (\cite{IS} Theorem 2.1   p 1079 ) \label{ThmIS}
{\begin{enumerate}
\item[$H1-$] $\mathcal{F}\in \mathcal{C}([0,T]\times \overline{\Omega}\times \R\times( \R^2\setminus\{0_{\R_2}\})\times \mathcal{S}_2)$,
\item[$H2-$] There exists $\gamma\in \R$  such that for every $(t,x,p,X)\in[0,T]\times \overline{\Omega}\times ( \R^2\setminus\{0_{\R_2}\})\times \mathcal{S}_2$,  the function  $\lambda \mapsto \mathcal{F}(t,x,\lambda, p,X)-\gamma \lambda $ is non-decreasing  on $\R$.
\item[$H3-$]  For every $R>0$,  there exists a continuous non-decreasing function  $\varphi_R:[0,+\infty[\to [0,+\infty[$such that  $\varphi_R(0)=0$ and  for every  $X,Y\in \mathcal{S}_2$ and $\mu_1, \mu_2\in [0,+\infty[$  satisfying : 
$$\left(\begin{array}{cc}X & 0 \\0 & Y\end{array}\right)\leq \mu_1\left(\begin{array}{cc}I & -I \\-I & I\end{array}\right)+\mu_2\left(\begin{array}{cc}I & 0 \\0 & I\end{array}\right),$$ then 
\begin{equation}
\begin{split}
& \mathcal{F}(t,x,\lambda,p,X)-\mathcal{F}(t,y,\lambda,q,-Y)  \geq \\& -\varphi_R\left(\mu_1(|x-y|^2+\rho(p,q)^2) +\mu_2+|p-q|+|x-y|(1+\max(|p|,|q|))\right),
\end{split}
\end{equation}
$\forall t\in[0,T]$, $x,y\in \overline{\Omega}$, $\lambda\in\R$ such that  $|\lambda| \leq R$ and $p,q \in\R^2\setminus\{0_{\R^2}\}$.
\item[$H4-$]  $\mathcal{F}$ is continuous at  $(t,x,\lambda,0,0)$ pour tout $(t,x,\lambda)\in[0,T]\times \overline{\Omega}\times \R$ in the following sense: $$-\infty < \mathcal{F}_*(t,x,\lambda,0,0)=\mathcal{F}^*(t,x,\lambda,0,0)<+\infty$$
where $\mathcal{F}^*$ (respectively  $\mathcal{F}_*$ ) are the upper (respectively lower) semi-continuous envelopes of  $\mathcal{F}$,  defined on $[0,T]\times \overline{\Omega}\times \R\times\R^2\times \mathcal{S}_2.$
\item[$B1-$] $B\in \mathcal{C}(\R^2\times \R^2)\cap \mathcal{C}^{1,1}(\R^2\times (\R^2\setminus\{0_{\R^2}\})).$
\item[$B2-$]  Pour tout $x\in \R^2$,  the function  $p\mapsto B(x,p)$  1-positively homogeneous with respect to  $p$, i.e., $B(x,\lambda p)=\lambda B(x,p), \forall \lambda\geq 0, p \in\R^2\setminus\{0_{\R^2}\}.$
\item[$B3-$]   There exists a  positive constant $\theta$ such that $\proscal{\nu(z)}{D_pB(z,p)}\geq \theta$ for every  $z\in \partial \Omega$ and  $p\in\R^2\setminus\{0_{\R^2}\}$. 
\\ Here  $\nu(z)$  is the unit outer normal vector of $\Omega$ at $z\in \partial \Omega$.

\end{enumerate}
Assume [$H1$, $H2$, $H3$, $H4$, $B1$, $B2$, $B3$] are satisfied and consider the following problem   
\begin{equation}(\mathcal{S})\label{eqn:thIS}
\left\{\begin{array}{ll}
\phi _t+\mathcal{F}(t,x,\phi,\nabla \phi,D^2\phi)=0, &\mbox{in }  \mathcal{Q}:=]0,T[ \times \Omega   \\
B(x,\nabla \phi)=0 &\mbox{on }  ]0,T[ \times \partial \Omega   \end{array} \right.
\end{equation}
 $\bullet$ Let $\phi \in USC ([0,T[\times \overline{\Omega})$  and $\psi \in LSC([0,T[\times \overline{\Omega})$ be, respectively, viscosity sub and supersolutions  of  (\ref{eqn:thIS}). 
If  $\phi(0,x)\leq \psi (0,x)$ for $x \in \overline{\Omega}$, then $\phi \leq  \psi $ sur $]0,T[\times \overline{\Omega}$.
 \\
 $\bullet $ For every function  $g\in\mathcal{C}(\overline{\Omega})$ there exists a unique viscosity solution  $\phi\in\mathcal{C}([0,T[\times \overline{\Omega})$  of   (\ref{eqn:thIS}) such that $ \phi(0,x)=g(x) $ on $\overline{\Omega} $.
 }} 
 
Here $USC([0,T[\times \overline{\Omega})$  (respectively $LSC([0,T[\times \overline{\Omega})$ denote the set of upper semicontinuous (respectively lower semicontinuous) functions.
\subsection{Existence and uniqueness of the solution of the evolution problem }\label{exist-uniq}

In the sequel we assume that the balloon force $H$ verifies $(\mathcal{H}_2)$ : $$(\mathcal{H}_2) : \qquad  H \mbox{ is Lipschitz continuous on } \overline{\Omega}.$$ 
In order to use   Ishii-Sato theorem  we have to  verify  every hypothesis. The Hamiltonian
  $\mathcal{F}(t,x,r,p, X)=F({\bf{x}},p,X)$  is 
defined by (\ref{eqn:F}): 

\begin{equation}\label{eqn:Ftilde}
F({\bf{x}},p,X) =\widetilde{g}({\bf{x}}) H({\bf{x}}) |p|
-\beta \widetilde{g}({\bf{x}})   \tr(A(p)X)
+\widetilde{g}({\bf{x}})  (1-|H({\bf{x}})|)\proscal{\hat{V}({\bf{x}})}{p}~.
\end{equation}
 Let us define the symmetric, semi-definite positive matrix $A({\bf{x}},p)$ : 
$$A({\bf{x}},p)=\beta\, \widetilde{g}({\bf{x}})(I- \frac{p\otimes p}{|p|^2}),$$
so  that  
\begin{equation}\label{eqn:H2}
F({\bf{x}},p,X)=-\tr(A({\bf{x}},p)X)+\widetilde{g}({\bf{x}}) H({\bf{x}}) |p|
 +\widetilde{g}({\bf{x}})(1-|H({\bf{x}})|)\proscal{\hat{V}({\bf{x}})}{p}~.
\end{equation}
In this case $F$ does not depend neither on $t$ nor   $\lambda$. 
\\
We choose a  Neumann-type boundary condition : 
 $$\proscal{\nabla \phi({\bf{x}})}{\nu({\bf{x}})}_{\R^2}=\displaystyle{\frac{\partial \Phi}{\partial \nu}({\bf{x}})}=0 \hbox{ on } \Gamma=\partial \Omega$$
  
\begin{enumerate}
\item [{[H1]}]  Function  $F \in \mathcal{C}(\overline{\Omega}\times( \R^2\setminus\{0_{\R_2}\})\times \mathcal{S}_2$.
 $p=0$ is a singular point.  \\ Theorem \ref{regulGVF} yields that the gradient vector flow $V$ is continuous on $[0,T] \times \overline{\Omega}$.
 
 In addition   we assumed that  the balloon function $H$ is continuous on $\overline{\Omega}$ .
 \\
 Therefore, the Hamiltonian $F$ is continuous on $\overline{\Omega}\times \R^2\setminus\{0_{\R^2}\} \times \mathcal{S}_2$
 
\item [{[H2]}]   Let us show there exists a constant $\gamma\in \R$ such that for each $(x,p,X)\in \times \overline{\Omega}\times ( \R^2\setminus\{0_{\R^2}\})\times \mathcal{S}_2$, the function $\lambda\mapsto F(x,p,X)-\gamma \lambda$ is non-decreasing on $\R$. Since $F$ does not explicitly depend on $\lambda$,  any negative constant $\gamma$  is suitable.

\item [{[H3]}]  As $F$ does not depend on $t$ and $r$, we have to find  a continuous increasing function  
 $\varphi :[0,+\infty[\to [0,+\infty[$ satisfying $\varphi(0)=0$ such that if $X,Y\in \mathcal{S}_2$ and $\mu_1, \mu_2\in [0,+\infty[$ satisfy : 
\begin{equation}\label{plus1}
\left(\begin{array}{cc}X & 0 \\0 & Y\end{array}\right)\leq \mu_1\left(\begin{array}{cc}I & -I \\-I & I\end{array}\right)+\mu_2\left(\begin{array}{cc}I & 0 \\0 & I\end{array}\right),\end{equation}
 then $\forall {\bf{x}},{\bf{y}} \in \overline{\Omega}$, and $p,q \in\R^2\setminus\{0_{\R^2}\}$
$$
 F({\bf{x}},p,X)-F({\bf{y}},q,-Y)  $$
 $$\geq  -\varphi \left(\mu_1(|{\bf{x}}-{\bf{y}}|^2+\rho(p,q)^2) +\mu_2+|p-q|+|{\bf{x}}-{\bf{y}}|(1+\max(|p|,|q|))\right), $$
 Let use   the following lemma \cite{LG} : 
{\lemma \label{ro}
If  $p,q\in\R^N\setminus\{0_{\R^2}\}$, then   : 
$$\left |\frac{p}{|p|}-\frac{q}{|q|}\right|\leq \frac{|p-q|}{\min(|p|,|q|)}:=\rho(p,q).$$}

Given $X,Y\in \mathcal{S}_2$ and  $\mu_1, \mu_2\in [0,+\infty[$ verifying  (\ref{plus1}).
 Let be $r,s\in\R^2$, so we have
$$\proscal{Xr}{r}+\proscal{Ys}{s}\leq \mu_1|r-s|^2+\mu_2(|r|^2+|s|^2)$$

 Let ${\bf{x}},{\bf{y}} \in \overline{\Omega}$ and  $p,q\in \R^2\setminus\{0_{\R^2}\}$. 
 Following   C. Le Guyader \cite{LG}, we split $F$ in three terms and verify [H3] for each term. 
\\
$ F({\bf{x}},p,X)-F({\bf{y}},q,-Y) = $
\begin{equation}\label{plus2} \begin{array}{ll}
&-\underbrace{\left(\tr(A({\bf{x}},p)X)+\tr(A({\bf{y}},q)Y)\right)}_{(a)} \\[0.3cm]
& + \underbrace{\widetilde{g}({\bf{x}})H({\bf{x}})|p|-\widetilde{g}({\bf{y}})H({\bf{y}})|q|}_{(b)}\\[0.3cm]
&  + \underbrace{\widetilde{g}({\bf{x}})(1-|H({\bf{x}})|)\proscal{\hat{V}({\bf{x}})}{p}_{\R^2} -\widetilde{g}({\bf{y}})(1-|H({\bf{y}})|)\proscal{\hat{V}({\bf{y}})}{q}_{\R^2}}_{(c)}~.\end{array}
\end{equation}
$\bullet$ \textbf{(a) estimate.- }
As  $A({\bf{x}},p)=\sigma({\bf{x}},p)\sigma^t({\bf{x}},p),$ 
$$
\tr(A({\bf{x}},p)X)+\tr(A({\bf{y}},q)Y) \leq  $$
$$ \mu_1 \tr\left((\sigma({\bf{x}},p)-\sigma({\bf{y}},q))(\sigma({\bf{x}},p)-\sigma({\bf{y}},q))^t \right)  +\mu_2  \left(\tr(\sigma({\bf{x}},p)\sigma^t({\bf{y}},p))+\tr(\sigma({\bf{y}},q)\sigma^t({\bf{y}},q))\right) $$
$$ \leq \mu_1 \tr \underbrace{\left((\sigma({\bf{x}},p)-\sigma({\bf{y}},q))(\sigma({\bf{x}},p)-\sigma({\bf{y}},q))^t \right)}_{(a1)}  +2\mu_2 \beta \delta$$ 
where $\delta$ is an upper bound of  $ \widetilde{g}$ on $\overline{\Omega}$ (for example $\displaystyle{\frac{1}{\sqrt{2\pi}\sigma} }$). Expression (a1)  verifies : 
\begin{equation*}
\begin{split}
&\tr \left((\sigma({\bf{x}},p)-\sigma({\bf{y}},q))(\sigma({\bf{x}},p)-\sigma({\bf{y}},q))^t \right)
\\& =\tr \left(\sigma({\bf{x}},p)\sigma^t({\bf{x}},p) -\sigma({\bf{x}},p)\sigma^t({\bf{y}},q)-\sigma({\bf{y}},q)\sigma^t({\bf{x}},p)+\sigma({\bf{y}},q)\sigma^t({\bf{y}},q)\right)
\\& =\tr \left(A({\bf{x}},p)-\sigma({\bf{x}},p)\sigma^t({\bf{y}},q)-\sigma({\bf{y}},q)\sigma^t({\bf{x}},p)+A({\bf{y}},q)\right)
\\&= \beta \widetilde{g}({\bf{x}})+2\beta\sqrt{\widetilde{g}({\bf{x}})}\sqrt{\widetilde{g}({\bf{y}})}\left(-2+\frac{1}{|p\|q|}(p_1q_1+p_2q_2)\right)+\beta \widetilde{g}({\bf{y}})
\end{split}
\end{equation*}
and with  $$\left(-2+\frac{1}{|p\|q|}(p_1q_1+p_2q_2)\right)\leq -\frac{1}{|p\|q|}(p_1q_1+p_2q_2)$$
we get : 
\begin{equation*}
\begin{split}
&\tr \left((\sigma({\bf{x}},p)-\sigma({\bf{y}},q))(\sigma({\bf{x}},p)-\sigma({\bf{y}},q))^t \right)
\\& \leq \beta \widetilde{g}({\bf{x}})-2 \beta \sqrt{\widetilde{g}({\bf{x}})}\sqrt{\widetilde{g}({\bf{y}})}\frac{1}{|p\|q|}(p_1q_1+p_2q_2)+\beta \widetilde{g}({\bf{y}})
\\&=\left|\sqrt{\beta \widetilde{g}({\bf{x}})}\frac{p}{|p|}-\sqrt{\beta \widetilde{g}({\bf{y}})}\frac{q}{|q|}\right|^2
\end{split}
\end{equation*}
We finally obtain 
\begin{equation} \tr(A({\bf{x}},p)X)+\tr(A({\bf{y}},q)Y) \leq \mu_1 \left|\sqrt{\beta \widetilde{g}({\bf{x}})}\frac{p}{|p|}-\sqrt{\beta \widetilde{g}({\bf{y}})}\frac{q}{|q|}\right|^2
 +2\mu_2 \beta \delta
\end{equation}
Moreover we have the relation : 
 $$ \begin{array}{ll}
  \displaystyle{\left|\sqrt{\beta \widetilde{g}({\bf{x}})}\frac{p}{|p|}-\sqrt{\beta \widetilde{g}({\bf{y}})}\frac{q}{|q|}\right|^2} 
  &\leq  
  \displaystyle{\left|\left(\sqrt{\beta \widetilde{g}({\bf{x}})}-\sqrt{\beta \widetilde{g}({\bf{y}})}\right)\frac{p}{|p|}+\sqrt{\beta \widetilde{g}({\bf{y}})}\left(\frac{p}{|p|}-\frac{q}{|q|}\right)\right|^2 }\\
 & \displaystyle{\leq 2\beta \left(\sqrt{\widetilde{g}({\bf{x}})}-\sqrt{\widetilde{g}({\bf{y}})}\right)^2+2\beta \widetilde{g}({\bf{y}})\left|\frac{p}{|p|}-\frac{q}{|q|}\right|^2.}\end{array}$$ 
%
So we get by lemma \ref{eqn:l2} : 
\begin{equation*}
\begin{split}
 \left|\sqrt{\beta \widetilde{g}({\bf{x}})}\frac{p}{|p|}-\sqrt{\beta \widetilde{g}({\bf{y}})}\frac{q}{|q|}\right|^2
&= \left|\sqrt{\beta \widetilde{g}({\bf{x}})}\frac{p}{|p|}-\sqrt{\beta \widetilde{g}({\bf{y}})}\frac{p}{|p|}+\sqrt{\beta \widetilde{g}({\bf{y}})}\left(\frac{p}{|p|}-\frac{q}{|q|}\right)\right|^2
\\& \leq 2\beta K_1^2|{\bf{x}}-{\bf{y}}|^2+2\beta \delta \left|\frac{p}{|p|}-\frac{q}{|q|}\right|^2.
\end{split}
\end{equation*}
 Thus we  conclude  with  lemma \ref{ro}.  
 
\begin{equation}\label{eqn:T2}
\tr(A({\bf{x}},p)X)+\tr(A({\bf{y}},q)Y)\leq \mu_1\left( 2\beta K_1^2 |{\bf{x}}-{\bf{y}}|^2+2\beta \delta 4 \rho(p,q)^2\right)+2\mu_2 \beta \delta 
\end{equation}
$\bullet$\textbf{ (b) estimate.- }
We have assumed   $H$ to be  Lipschitz continuous on $\overline{\Omega}$ so  ${\bf{x}}\in \overline{\Omega} \mapsto \widetilde{g}({\bf{x}})H({\bf{x}})$ is Lipschitz continuous and bounded as well. 
 Let be ${\bf{x}},{\bf{y}}\in \overline{\Omega}$
$$ \begin{array}{ll}
 \left| \widetilde{g}({\bf{x}})H({\bf{x}})|p|-\widetilde{g}({\bf{y}})H({\bf{y}})|q|\right|
 & \leq \left|\left(\widetilde{g}({\bf{x}})H({\bf{x}})-\widetilde{g}({\bf{y}})H({\bf{y}})\right)|p|\right|
  + \left|\widetilde{g}({\bf{y}})H({\bf{y}})\|p|-|q\|\right| \\
 & \leq K_2|{\bf{x}}-{\bf{y}}|\max(|p|,|q|)+\theta \|p|-|q\|
\\& \leq K_2|{\bf{x}}-{\bf{y}}|\max(|p|,|q|)+\theta |p-q|
\end{array}
$$
where
$K_2$ is the Lipschitz-constant of  the function ${\bf{x}}\in \overline{\Omega} \mapsto \widetilde{g}({\bf{x}})H({\bf{x}})$ and
$\theta$ is a bound of ${\bf{x}}\in \overline{\Omega}  \mapsto  |\widetilde{g}({\bf{x}})H({\bf{x}})|$.

$\bullet$\textbf{ (c) estimate.- }  
We have assumed $(\mathcal{H}_2)$ so ${\bf{x}}\in \overline{\Omega} \mapsto \widetilde{g}({\bf{x}})(1-|H({\bf{x}})|)\hat{V}({\bf{x}})$ is Lipschitz continuous and bounded.

\begin{equation*}
\begin{split}
&\left| \widetilde{g}({\bf{x}})(1-|H({\bf{x}})|)\proscal{\hat{V}({\bf{x}})}{p}_{\R^2}  -\widetilde{g}({\bf{y}})(1-|H({\bf{y}})|)\proscal{\hat{V}({\bf{y}})}{q}_{\R^2} \right|
\\& \leq  \left|\proscal{\widetilde{g}({\bf{x}})(1-|H({\bf{x}})|)\hat{V}({\bf{x}})-\widetilde{g}({\bf{y}})(1-|H({\bf{y}})|) \hat{V}({\bf{y}})}{p}_{\R^2} \right|\\
&\qquad \qquad +\left|\proscal{\widetilde{g}({\bf{y}})(1-|H({\bf{y}})|) \hat{V}({\bf{y}})}{p-q}_{\R^2}   \right|
\\& \leq K_3 |{\bf{x}}-{\bf{y}}|\max(|p|,|q|)+\zeta|p-q|
\end{split}
\end{equation*}

where
$K_3$ is the Lipschitz-constant of  the function ${\bf{x}}\in \overline{\Omega} \mapsto \widetilde{g}({\bf{x}})(1-|H({\bf{x}})|)\hat{V}({\bf{x}})$ and
$\zeta$ is a bound of ${\bf{x}}\in \overline{\Omega}  \mapsto |\widetilde{g}({\bf{x}})(1-|H({\bf{x}})|)\hat{V}({\bf{x}})|.$
 \\
Finally  relation (\ref{plus2}) gives 
$$ \begin{array}{ll}
 -\left(F({\bf{x}},p,X)-F({\bf{y}},q,-Y)\right) &\leq \mu_1\left( 2\beta K_1^2|{\bf{x}}-{\bf{y}}|^2+8\beta \delta \rho(p,q)^2\right)+2\mu_2 \beta \delta
\\& +(K_2|{\bf{x}}-{\bf{y}}|\max(|p|,|q|)+\theta |p-q|)
\\& +(K_3 |{\bf{x}}-{\bf{y}}|\max(|p|,|q|)+\zeta|p-q|).
\\[0.5cm]
&\leq \max(2 \beta K_1^2,8 \beta \delta,K_2+K_3,\theta+\zeta)
[(\mu_1( |{\bf{x}}-{\bf{y}}|^2+\rho(p,q)^2) \\ &
+\mu_2) +\max(|p|,|q|)|{\bf{x}}-{\bf{y}}|
+|p-q|]  \\[0.5cm]
&\leq  \max(2 \beta K_1^2,8 \beta \delta,K_2+K_3,\theta+\zeta)
 \left(\mu_1( |{\bf{x}}-{\bf{y}}|^2+\rho(p,q)^2)\right. \\ &
+\left. \mu_2+(1+\max(|p|,|q|))|{\bf{x}}-{\bf{y}}|
+|p-q| \right)~.
\end{array} $$
Finally
 \begin{equation*}
\begin{split}
&F({\bf{x}},p,X)-F({\bf{y}},q,-Y)\\&\geq -\varphi \left(\mu_1( |{\bf{x}}-{\bf{y}}|^2+\rho(p,q)^2)
+\mu_2+(1+\max(|p|,|q|))|{\bf{x}}-{\bf{y}}|
+|p-q| )\right)
\end{split}
\end{equation*}
  where the function $\varphi$ is defined by : $$\varphi(m)= \max(2 \beta K_1^2,8 \beta \delta,K_2+K_3,\theta+\zeta) m$$
Hypothesis $[H3]$ is then verified.
 
\item [{[H4]}]  
$F$ is continuous at $({\bf{x}},0,0)$ for any ${\bf{x}} \in \overline{\Omega}$ because
$ F_*({\bf{x}},0,0)=F^*({\bf{x}},0,0)=0$. 

\vspace{0.7cm}

\nd For the three last hypothesis $B1$, $B2$ and $B3$  on the boundary condition, the proof is the same as in  C. Le Guyader \cite{LG}.
\item  [{[B1]}] This hypothesis consists in showing that the function $B$ which defines the boundary condition on $]0,+\infty[ \times \partial \Omega $ is $ \mathcal{C}(\R^N\times \R^N)\cap \mathcal{C}^{1,1}(\R^N\times (\R^N\setminus\{0_{\R^N}\})).$ We have chosen a Neumann-type condition, by denoting $\nu({\bf{x}})$ the outward unit normal vector of $\partial \Omega$ at the point ${\bf{x}}$, our boundary condition is written :  $$B({\bf{x}},p)=\proscal{\nu({\bf{x}})}{p}_{\R^2}.$$
The point $[B1]$ will be satisfied if $\nu$ is a $\mathcal{C}^{1,1}$ vector  field which is the case if $\Omega $ has a $\mathcal{C}^2$ boundary. 
\item [{[B2]}]  $B$ is  1-positively homogeneous with respect to  $p$ : $$B({\bf{x}},\lambda p)=\proscal{\nu({\bf{x}})}{\lambda p}_{\R^2}=\lambda B({\bf{x}},p),\forall \lambda \geq 0, p\in{\R^2\setminus\{0_{\R^2}\}}.$$
\item [{[B3]}]  Let $z\in\partial \Omega$.$$\proscal{\nu(z)}{D_pB(z,p)}_{\R^2}=|\nu(z)|^2=1.$$
The last condition is verified with $\theta=1$, 
\end{enumerate}
We may now conclude since assumptions [ $H1$, $H2$, $H3$, $H4$, $B1$, $B2$ , $B3$]  are satisfied. 
{\theoreme{Assume that the image  function $I$ verifies the hypothesis $(\mathcal{H}_I)$ and the ballon force  $H$ satisfies $(\mathcal{H}_2)$. 
Consider the following problem 
\begin{equation}\label{bis}
\hspace{-0.5cm}\left\{\begin{array}{l}
 \Phi_t(t,{\bf{x}})-g(|\nabla I({\bf{x}})|^2) \left( (\beta \kappa ({\bf{x}})-H({\bf{x}}))|\nabla \Phi ({\bf{x}})| - (1- |H({\bf{x}})|)\proscal{ \hat V ({\bf{x}})}{\nabla \Phi({\bf{x}})} \right) =0 \\
 \hfill\mbox{in }  ]0,+\infty[ \times \Omega , \\
 \displaystyle{\frac{\partial \Phi}{\partial \nu}({\bf{x}})=0} ~\mbox{ on  } ~~ ]0,+\infty[ \times\partial \Omega~,
\end{array} \right.
\end{equation}
$\bullet $ 
  Let $\Phi \in USC ([0,T[\times \overline{\Omega})$ and $\Psi \in LSC ([0,T[\times \overline{\Omega})$ be, respectively, viscosity sub and supersolutions of :  
If $\Phi(0,x)\leq \Psi (0,x)$ for $x\in \overline{\Omega}$, 
then $\Phi \leq  \Psi $ in $]0,T[\times \overline{\Omega}$.
 \\
 $\bullet$ For  every $g\in\mathcal{C}(\overline{\Omega})$, there is a unique viscosity solution $\Phi\in\mathcal{C}([0,T[\times \overline{\Omega})$ of (\ref{bis})  satisfying   $\Phi(0,x)=g(x) $ on $\Omega$.\\
 $\bullet$ Equation (\ref{eqn:levelsetparagios})  has  a unique viscosity solution $\Phi \in \mathcal{C}([0,T[\times \overline{\Omega})$. }}
 
\section{Conclusion}
We have  recalled   Gradient Vector Flow model that we hahe justified and we have  given regularity properties.  Then we proved existence and uniqueness fo viscosity solutions of the Hamilton-Jacobi equation derived fron the GVF-geodesic active contour model. 

Next step is to perform the numerical realization of this GVF-geodesic active contour process solving   he Hamilton-Jacobi equation (\ref{bis}). We shall use his method to the perform  ``tuffeau'' tomographic images segmentation. This material is has been used during past centuries to build monuments as castles and churches  in the  Val-de-Loire area.  These images allow to get information on the structure of the damaged material : we have to identify different phases as calcite (light grey), silice (dark grey) and porosity (black). 

We shall combine the GVF-geodesic model with a region segmentaion  approach   to identify the three constituents of tuffeau . The segmentation of these images is a step of pretreatment which aims at reconstructing the stone porosity domain.


\end{document}